\documentclass[11pt]{article}
\usepackage{amsfonts}
\usepackage{amsmath}
\usepackage{epsfig}
\usepackage[compress]{cite}
\usepackage{graphics}
  
\begin{document}

\begin{center}
\bf{A random map implementation of implicit filters} 
\end{center}

\begin{center}
 Matthias Morzfeld$^{1}$, Xuemin Tu$^2$, Ethan Atkins$^{1,3}$, Alexandre J. Chorin$^{1,3}$\\
\end{center}
\vspace{5mm}
1. Lawrence Berkeley National Laboratory, Berkeley, CA \\
2. Department of Mathematics, University of Kansas, Lawrence, KS \\ 
3. Department of Mathematics, University of California, Berkeley, CA \\
\vspace{5mm}

\begin{center}
\bf{Abstract}
\end{center}

Implicit particle filters for data assimilation generate high-probability samples by representing each particle location as a separate function of a common reference variable. This representation requires that a certain underdetermined equation be solved for each particle and at each time an observation becomes available. We present a new implementation of implicit filters in which we find the solution of the equation via a random map. As examples, we assimilate data for a stochastically driven Lorenz system with sparse observations and for a stochastic Kuramoto-Sivashinski equation with observations that are sparse in both space and time.

\vspace{5mm}
\noindent

Keywords: data assimilation; implicit sampling; particle filters; sequential Monte Carlo\\

AMS Subject Classification: 60G35, 62M20, 86A05

\section{Introduction}
\label{sec:Introduction}
In many applications in science and engineering the state of a system must be identified from an uncertain model supplemented by a stream of noisy data. Problems of this kind are typically formulated in terms of an It\^{o} stochastic differential equation (SDE)
\begin{equation}
	\label{eq:SODE}
	dx(t)=f(x,t)dt+g(x,t)dW,
\end{equation}
where $t\geq 0$ is the independent variable, the state $x(t)$ is a real $m$-dimensional column vector, $f(x,t)$ is a real $m$-dimensional vector function $g(x,t)$ is a real $m\times m$ matrix and $dW$ is an $m$-dimensional Brownian motion (BM). The probability density function (pdf) of the state at $t=0$ is known. As the solution of the SDE evolves, measurements 
\begin{equation}
	\label{eq:Observation}
	b^n = h(x^n,t^n)+QV^n
\end{equation}
are recorded at times $t^n, n=1,2,3,\ldots$, where $h$ is a $k$-dimensional vector function ($k\leq m$), $Q$ is a real $k\times k$ matrix, and $V^n$ is a $k$-dimensional vector, whose components are $k$ independent standard normal variates.  We assume for the remainder of this paper that $V^n$ is independent of the BM in (\ref{eq:SODE}). The goal is to use both the model (\ref{eq:SODE}) and the observations (\ref{eq:Observation}) to determine the state of the system.

The best estimate of the state $x^n=x(t^n)$ is, under wide conditions,  the mean of the probability density defined by the SDE and conditioned on the data.
In practice, the SDE must be discretized,
as described, for example, in~\cite{Doucet2001}, so that we are dealing with a discrete
recursion conditioned by discrete data.   
If the model (\ref{eq:SODE}) as well as the observations (\ref{eq:Observation}) are linear and if, in addition, the initial data are Gaussian, the conditional expectation can be computed via the Kalman-Bucy filter \cite{Kalman1961}. Strategies for tackling nonlinear, non-Gaussian problems include the ensemble Kalman filter \cite{Evensen}, the extended Kalman filter \cite{Gelb1974,Jazwinski1970}, the unscented Kalman filter \cite{Julier1997}, and variational methods \cite{ProvostSalmon1986,Sasaki1958,Sasaki1970}. These data assimilation strategies require Gaussian approximations or a linearization of the model and sometimes yield rather poor results if the nonlinearity is strong. We refer to \cite{Ide1997,Jardak2010,Miller1994,Bocquet2010,VanLeeuwen2009} for a review of various data assimilation algorithms, their applications and limitations.  

Particle filters \cite{Doucet2002,Doucet2001,Doucet2000, Gordon1993,VanLeeuwen2010,Weare2009} are sequential Monte Carlo tools that do not rely upon Gaussianity or linearity assumptions. In particle filters one works with a collection of ``particles" (replicas of the system), whose empirical distribution approximates the conditional pdf at the $n$-th step. One moves all particles forward in time using some guess of the pdf at the next step, and then one corrects the guess by weighting the particles.  The procedure is repeated at the next time an observation becomes available. The catch is that it is difficult to guess the next density 
accurately before doing the calculations; with most weighting schemes, many of the weights are therefore very small so that most of the computational effort is wasted on unlikely particles. As a consequence, the number of particles required can grow catastrophically with the dimension of the SDE~\cite{Bickel2008, Snyder2008}. 

The implicit filter \cite{Chorin2010b,Chorin2009} is designed to remedy this problem, i.e. it attempts to make nonlinear data assimilation feasible in high dimensional SDEs. The basic idea is to reverse the standard procedure. Rather than generating a sample and then computing its probability, the implicit filter finds regions of high probability taking the observations into account, and then looks for samples in these regions; this generates a thin beam of high probability particles, focussed on the observations, making the number of required particles manageable. 
The focusing is done by connecting the samples to a fixed reference density through a map
that satisfies a data-dependet algebraic equation. This map is not unique and the efficiency of the sampling depends on the map one chooses. In the present paper we present an efficient implementation in which the reference variables are connected to the samples by a random map. We demonstrate its effectiveness on two test problems: a stochastically driven Lorenz model with sparse data and a stochastic Kuramoto-Sivashinski equation with data sparse in both space and time. We compare the implicit filter to a Sampling-Importance-Resampling (SIR) filter (see \cite{Doucet2002,Doucet2001,Gordon1993}) which constructs a prior density using the SDE and only later re-weights the particle positions by the observations. 

A different algorithm for guiding particles towards high probability regions has been presented in \cite{VanLeeuwen2010}. If the observations are linear and available at every time step, our filter reduces to a variant of the optimal SIR filter \cite{Doucet2002,Bocquet2010,DelMoral1998}, as described below.
 
\section{Implicit sampling: basic ideas}
\label{sec:ImplicitSampling}
We start by reviewing the general framework of implicit sampling (see 
 \cite{Chorin2010b,Chorin2009}), explicitly allowing for the possibility that the observations are sparse in time, i.e., not necessarily available at every time step.
We assume that the 
SDE~(\ref{eq:SODE}) has been approximated by a difference scheme with time step $\delta$, in the form:
\begin{equation}
	\label{eq:SODEDiscrete}
	x^{n+1}=R(x^n,t^n)+G(x^n,t^n)\Delta W^{n+1},
\end{equation}
where the functions $R$ and $G$ depend on the scheme we use,  $\Delta W^{n+1}\sim N(0,\delta)$ and $t^n$ is shorthand for $n\delta$. For more details on the discretization, see \cite{KloedenPlaten} and the examples below. For simplicity, we assume a constant time step $\delta$. The generalization of implicit sampling to higher order integration schemes is straightforward, see \cite{Chorin2010b} and below.

Assume we are given a collection of $M$ particles with positions $X_j^n$, $j=1,\dots,M$, whose empirical density approximates the conditional pdf at time $t^n$, and suppose the next observation is available after $r$ time steps, $r$ a positive integer. Bayes' theorem can be used to show that the pdf of the \emph{j}th particle at times $t^0,\dots,t^n,t^{n+1},\dots,t^{n+r}$, conditioned on the available observations $b^{1},\dots,b^{n+r}$, is 
\begin{eqnarray}
	\label{eq:JointPDF}
	P(X_j^{0,\dots,n+r}\mid b^{1\dots n+r})&=& P(X_j^{0,\dots,n} \mid b^{1,\dots,n} )\nonumber \\ 
	 & &\times P(X_j^{n+1} \mid X_j^{n} )  \ldots   P(X_j^{n+r-1} \mid X_j^{n+r-2})\nonumber\\
	  & &\times P(X_j^{n+r}\mid X_j^{n+r-1}) P(b^{n+r}\mid X_j^{n+r})/Z^n
\end{eqnarray}
where
$Z^n$ is a normalization constant independent of the particles and $X_j^{0,\dots,n}$ is an abbreviation for $X_j^0,X_j^1,\dots,X_j^n$. Implicit sampling is a recipe for obtaining high probability samples from (\ref{eq:JointPDF}). 

For ease of notation we introduce the shorthand notation $X_j$ for the $r\cdot m$ dimensional column vector $[(X_j^{n+1})^T,\dots,(X_j^{n+r})^T]^T$ (the state trajectory of the particle) and define, for each particle, a function $F_j(X_j)$ by
\begin{eqnarray}
\exp(-F_j(X_j))&=& P(X_j^{n+1} \mid X_j^{n} ) \ldots  P(X_j^{n+r-1} \mid X_j^{n+r-2})\nonumber\\
\label{eq:F}
          &  & \times P(X_j^{n+r}\mid X_j^{n+r-1}) P(b^{n+r}\mid X_j^{n+r}).
   \end{eqnarray}
To obtain a sample we solve the algebraic equation
\begin{equation}
	\label{eq:UnderdeterminedEq}
	F_j(X_j)-\phi_j = \frac{1}{2}\xi_j^T\xi_j,
\end{equation}
where $\xi_j$ is a realization of a random variable $\xi$, drawn from a given, fixed reference density, say a $r\cdot m$ dimensional, multivariate normal distribution $P_{\xi}=\mathcal{N}(0,I)$. The additive, deterministic factor $\phi_j$ is needed to make the equation $F_j(X_j) = \xi_j^T\xi_j/2$ solvable (the left-hand-side is real, but the right-hand-side is non-negative). The choice 
\begin{equation}
	\label{eq:TheMin}
	\phi_j= \min F_j,
\end{equation}
where $\min F_j$ is the global minimum of $F_j$, will do the job, and this is the choice we make in the present paper. We solve (\ref{eq:UnderdeterminedEq}) for each particle because the functions $F_j$ vary from particle to particle due to different parameters~$X_j^{n}$. 

The variable $\xi$ on the right hand side of (\ref{eq:UnderdeterminedEq}) is known and easy to sample, and by definition most of its samples will be high-probability samples near the origin. By equations (\ref{eq:UnderdeterminedEq}) and (\ref{eq:TheMin}) the corresponding values of $F_j(X_j)$ will be near the minimum of $F_j$ and therefore will have a high probability, so that with high probability we will have high probability samples. The probability density of the samples $X_j$ corresponds to the ``prior density'' in the usual Bayesian
sampling, but it is not a prior density in the usual sense, because the new positions of the particles
are obtained by solving different equations, rather than by sampling a common prior. The prior here is a parametrized family of functions of the reference variable. 

The empirical density defined by the new particle positions differs from the target density so that each sample must be weighted by the ratio of its probability with respect to the target density to its proposal probability \cite{VanLeeuwen2009,VanLeeuwen2010,Doucet2002}. Using $P^{n+1}(X_j)$ and $P_{X_j}(X_j)$ as a shorthand notation for the target density (\ref{eq:JointPDF}) and the density defined by (\ref{eq:UnderdeterminedEq}) respectively, we can obtain the weight $w_j^{n+1}$ of the particle $X_j$ at time $t^{n+1}$, i.e. its probability with respect to the target density:
\begin{eqnarray}
	w_j^{n+1} & \propto & \quad \frac{P^{n+1}(X_j)}{P_{X_j}(X_j)}, \nonumber \\ 
	    & \propto &w_j^{n}\quad \frac{\exp(-F_j(X_j))}{P_\xi(\xi_j)} \mbox{  }J, \nonumber \\
	    &\propto &w_j^{n}\quad\frac{\exp(-0.5\xi_j^T\xi_j-\phi_j)}{\exp(-0.5\xi_j^T\xi_j)} \mbox{  } J, \nonumber \\
	    \label{eq:WeightGeneral}
	    & \propto &w_j^{n}\quad \exp(-\phi_j) \mbox{  }J,
\end{eqnarray}
where $J =  \left|\det \partial X_j/\partial \xi\right|$ is the Jacobian of the map. Having obtained the weights of all the particles, we normalize the weights so that their sum equals one. The variability of the weights modifies the reference density. The weights can be eliminated by resampling. Various resampling strategies and algorithms are discussed in \cite{Doucet2002}. Upon resampling, all particles have equal weights, so that, in particular, it is legitimate to omit the factor $w_j^{n}$ in equation (\ref{eq:WeightGeneral}).

For future reference, we rewrite the function $F_j$ in a slightly different form:
\begin{eqnarray}
F_j(X_j)&=& \frac{1}{2}\left(X_j^{n+1} -f(X_j^{n})\right)^T \left(G(X_j^n)G(X_j^n)^T\right)^{-1}\left(X_j^{n+1} -f(X_j^{n})\right)\nonumber\\
&+& \frac{1}{2}\left(X_j^{n+2} -f(X_j^{n+1})\right)^T \left(G(X_j^{n+1})G(X_j^{n+1})^T\right)^{-1}\left(X_j^{n+2} -f(X_j^{n+1})\right)\nonumber\\
& \vdots &  \nonumber\\
&+&\frac{1}{2} \left(X_j^{n+r} -f(X_j^{n+r-1})\right)^T \left(G(X_j^{n+r-1})G(X_j^{n+r-1})^T\right)^{-1}\left(X_j^{n+r} -f(X_j^{n+r-1})\right)\nonumber\\
          &  + & \frac{1}{2}\left(h(X_j^{n+r}) -b^{n+r})\right)^T \left(QQ^T\right)^{-1}\left(h(X_j^{n+r}) -b^{n+r})\right) + Z_j,
\end{eqnarray}
where $Z_j$ is a positive constant that can be computed from the normalization constants of the pdf's in the definition of $F_j$ in (\ref{eq:F}). This constant need not be computed because it drops out in (\ref{eq:UnderdeterminedEq}) when $\phi_j=\min F_j$.

Our construction can be readily generalized to SDE integration schemes with intermediate random steps. 
For example, suppose one is integrating the SDE (\ref{eq:SODE}) with additive noise, i.e. $g(x^n,t^n) = g = \mbox{constant}$,
using the Klauder-Petersen scheme \cite{Kla1}:
\begin{eqnarray}
\label{eq:KPScheme1}
x^{n+1,*}&=&x^n+\delta f(x^n)+g\sqrt{\delta}\Delta W_1,\\
\label{eq:KPScheme2}
x^{n+1}&=&x^n+\frac{\delta}{2}\left(f(x^n)+f(x^{n+1,*})\right)+g\sqrt{\delta}\Delta W_2,
\end{eqnarray}
where $\Delta W_1,\Delta W_2$ are $m$ dimensional Gaussians with mean zero and variance $I$. For simplicity, assume that observations
$b^n=h(x^n)+\eta_3$, with $\eta_3 \sim \mathcal{N}(0,s)$, are
available at every step. 
Dropping the index of the particles, the probability of the pair $(X^{n+1,*},X^{n+1})$ is proportional to
$\exp(-F)$, with
\begin{eqnarray}
	F&=& \frac{\left| \right|X^{n+1,*}-X^n-\delta f(X^n)\left|\right|^2}{2\delta g^2}
	 +\frac{\left| \right|X^{n+1}-X^n -\frac{\delta}{2}(f(X^n)+f(X^{n+1,*})\left| \right|^2}{2\delta g^2} \nonumber \\
	 \label{eq:KPF} 		 
	 & &+\frac{\left| \right| h(X^{n+1})-b\left| \right|^2}{2s}+Z, 
\end{eqnarray}
where the norm $||x||=\sqrt{x^Tx}$ is the Euclidean norm. All one has to do then is solve (\ref{eq:UnderdeterminedEq}) for the pair $(X^{n+1,*},X^{n+1})$,  with a sample $\xi_j$, drawn from a $2m$ dimensional Gaussian reference density, on the right-hand-side. 

The effectiveness of the filter rests on one's ability to solve the basic equation (\ref{eq:UnderdeterminedEq}) efficiently. This equation is underdetermined - it is a single equation connecting the $m\cdot r$ components of $X_j$ to the
reference variable $\xi$. Each solution algorithm defines a map from $\xi$ to $X_j$, and one has a great deal of freedom in choosing this map. Effective algorithms take advantage of this freedom. The conditions that the map must satisfy were derived and explained in \cite{Chorin2010b}: the map should be (\emph{i}) one-to-one and onto with probability one (so that the whole sample space is covered); (\emph{ii}) smooth near the high-probability region of $\xi$ (so that the weights do not vary unduly from particle to particle); (iii) it should map the neighborhood of zero onto a neighborhood of the minimum of $F$,  and (\emph{iv}) there should be an easy way to evaluate the Jacobian $\left|\det \partial X / \partial \xi\right|$. In our experience, condition (\emph{iv}) is often the most onerous to satisfy in nonlinear problems.

\section{Solution of the implicit sampling equation via a random map}
\label{sec:Construction}
A solution algorithm for equation (\ref{eq:UnderdeterminedEq}) defines a map from 
$\xi$ to the sample $X_j$. This map is not unique and should satisfy conditions (\emph{i}-\emph{iv}) above. Various ways to solve  (\ref{eq:UnderdeterminedEq}) have been presented in \cite{Chorin2010b,Chorin2009}.  In this section, we present a map that is random. 

First, we need to find the additive factor $\phi_j=\min F_j$ in (\ref{eq:UnderdeterminedEq}). We propose to find $\min F_j$ using standard tools, e.~g. Newton's method. It is also important to note that the Hessian of $F_j$ typically has a sparse block structure. This sparsity depends on the integration scheme we use for the discretization of the underlying SDE and can be exploited in the implementation of the algorithm. In the examples in sections \ref{sec:Lorenz} and \ref{sec:KS}, we present strategies for obtaining a ``good'' initialization for the minimization and then use a few straightforward Newton steps to polish the initial guess. We had no problems with this approach in the examples we considered. However, a quasi-Newton method can be used if the Hessian of $F_j$ is out of reach. More sophisticated minimization strategies, e.g. a trust-region method, may be preferable in other applications. 

In the present paper we assume that $F_j$ is convex and the Hessian evaluated at the global minimum $\phi_j$ is nonsingular. In \cite{Chorin2010b} we discussed what to do when this is not the case, in particular we presented strategies for replacing $F$ by a convex function without bias or loss. The methods presented here are compatible with the constructions in \cite{Chorin2010b}.

To find a sample $X_j$, we solve (\ref{eq:UnderdeterminedEq}) via the random Ansatz:
\begin{equation}
	\label{eq:Ansatz}
	X_j = \mu_j + \lambda_j  L_j^T\eta_j,
\end{equation}
where $\eta_j = \xi_j/\sqrt{\xi_j^T \xi_j}$, $\xi_j$ is a sample of the Gaussian reference variable $\xi\sim \mathcal{N}(0,I)$ and $\mu_j$ is the location of the minimum of $F_j$, i.e. $\mu_j=\mbox{argmin }F_j$. The invertible $rm\times rm$ matrix $L_j$ is deterministic, under our control, and remains to be chosen (see below). By substitution of (\ref{eq:Ansatz}) into (\ref{eq:UnderdeterminedEq}) we obtain a single algebraic equation in a single variable $\lambda_j$. The equation can be readily solved  and its solution defines the sample $X_j$. A data assimilation problem of arbitrary dimension thus boils down to a minimization of a known function followed by the solution of an algebraic equation in one variable. When $\lambda_j(\xi) \neq 0$ is continuous, the map (\ref{eq:Ansatz}) is one-to-one and onto almost surely so that requirement (\emph{i}) is satisfied.

\begin{figure}[h!t]
	\begin{center}
		{\includegraphics[width=.7\textwidth]{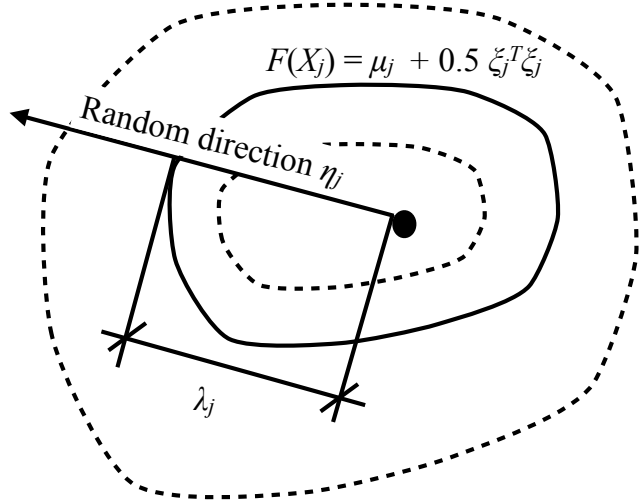}}
	\caption{Geometry of the random map $X_j\rightarrow \xi$.}
	\label{fig:Geometry}
	\end{center}
\end{figure}

The process of finding a sample via the random map (\ref{eq:Ansatz}) can be interpreted geometrically. Assuming that the level sets of $F_j$ are closed, the algebraic equation (\ref{eq:UnderdeterminedEq}) has a solution in every direction. We generate a random direction by sampling the reference density and computing $\eta$, which is uniformly distributed on the unit sphere. We determine how far we need to walk along the random direction $\eta$ to hit the level set $F(X_j) = \phi_j + 0.5\xi_j^T\xi_j$ by solving (\ref{eq:UnderdeterminedEq}) with the map~(\ref{eq:Ansatz}). The matrix $L_j$ is used to incorporate prior information. The geometry of the map is illustrated in figure \ref{fig:Geometry}. 

What remains to be done is compute the Jacobian of the map. The easiest way to do this calculation is as follows. We first rewrite~(\ref{eq:Ansatz}) as
\begin{equation}
	\label{eq:AnsatzHat}
	X_j = \mu_j + \hat{\lambda}_j  L_j^T\xi_j,
\end{equation}
with $\hat{\lambda}_j = \lambda_j/\sqrt{\xi_j^T\xi_j}$. From (\ref{eq:AnsatzHat}), we compute the derivatives
\begin{equation}
	\label{eq:Step1}
	\frac{\partial X}{\partial \xi}= L^T\left(\xi\frac{\partial \hat{\lambda}}{\partial \xi}\right) +\hat{\lambda} L^T,
\end{equation}
where we droped the index $j$ for the particles for convenience and where $\partial{\hat{\lambda}}/\partial{\xi}$ denotes the gradient (a row vector) of the scalar $\hat{\lambda}$ with respect to the reference variable $\xi$ (a column vector). By the chain rule, we obtain
\begin{equation}
	\label{eq:Step2}
	\frac{\partial \hat{\lambda}}{\partial \xi} = \frac{\partial \hat{\lambda}}{\partial \rho} \frac{\partial \rho}{\partial \xi} = 2\frac{\partial \hat{\lambda}}{\partial \rho} \xi^T,
\end{equation}
where $\rho = \xi^T\xi$, and substitute the result into (\ref{eq:Step1}) to get
\begin{equation}
\label{eq:Step3}
	\frac{\partial X}{\partial \xi}= L^T \left(2\frac{\partial \hat{\lambda}}{\partial \rho}\xi \xi^T+\hat{\lambda} I\right).
\end{equation}
Using standard rules for determinants such as $\det (A+BC) =\det A \cdot \det(I+CA^{-1}B)$, we calculate the Jacobian:
\begin{equation}
	\label{eq:JacobianInt}
	J =  |\det L| \mbox{ }  \left|\hat\lambda^{rm-1} \left(\hat\lambda+2\frac{\partial \hat\lambda}{\partial \rho}\rho \right)\right|.
\end{equation}
Substitution of $\lambda$ for $\hat\lambda$ finally gives
\begin{equation}
	\label{eq:Jacobian}
	J = 2  |\det L| \mbox{ } \rho^{1-rm/2} \mbox{ } \left|\lambda^{rm-1} \frac{\partial \lambda}{\partial \rho}\right|.
\end{equation}
A formula for the scalar derivative $\partial\lambda/\partial\rho$ can be obtained by implicit differentiation of (\ref{eq:UnderdeterminedEq}) combined with (\ref{eq:Ansatz})
\begin{equation}
	\frac{\partial \lambda_j}{\partial \rho_j}= \frac{1}{2\left(\nabla F_j\right) L_j^T\eta_j},
\end{equation}
where $\nabla  F_j$ denotes the gradient of $F_j$ (an $rm$-dimensional row vector). Alternatively, $d\lambda_j/d\rho_j$ can be computed numerically by putting $\lambda= \lambda+\Delta\lambda$, computing a new $X_j+\Delta X_j$ using Eq. (\ref{eq:Ansatz}), followed by evaluation of the left hand side of Eq. (\ref{eq:UnderdeterminedEq}) to get $\rho+\Delta\rho$, and differencing. The Jacobian (\ref{eq:Jacobian}) can thus be evaluated readily and condition (\emph{iv}) in section \ref{sec:ImplicitSampling} is satisfied. From (\ref{eq:WeightGeneral}), we compute the weights attached to each particle:
\begin{equation}
    \label{eq:WeightMap}
	w_j^{n+1} \propto w_j^{n}\exp(-\phi_j) |\det L_j| \mbox{ } \rho_j^{1-rm/2} \mbox{ } \left|\lambda_j^{rm-1} \frac{\partial \lambda_j}{\partial \rho_j}\right|.
\end{equation}

\begin{figure}[h!t]
	\begin{center}
		{\includegraphics[width=0.8\textwidth]{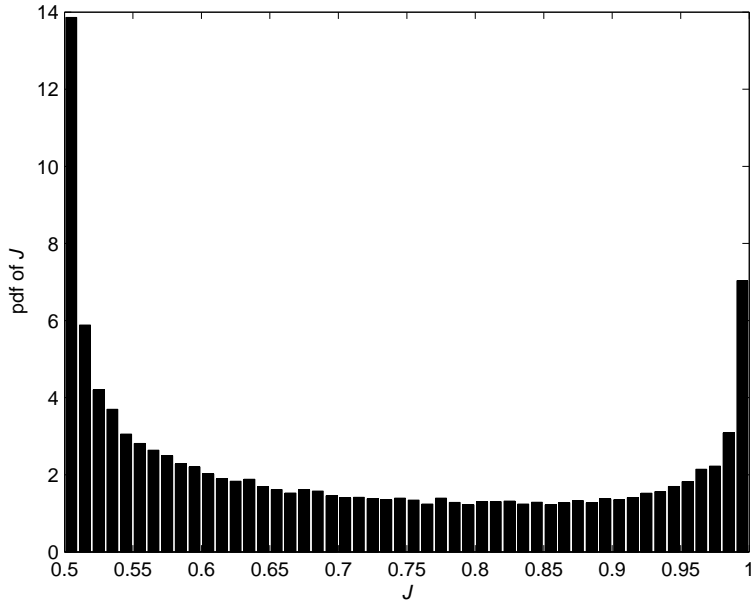}}
	\caption{Probability density function of the Jacobian for $n=2$ and a given~$\Sigma$.}
	\label{fig:BetaPDF}
	\end{center}
\end{figure}

We now need to choose the matrix $L_j$. In the examples we considered (see sections 4 and 5), the filters performed poorly  with the naive choice $L_j=I$. To understand why, suppose that observations are linear and available at every step. Equation (\ref{eq:Observation}) becomes 
\begin{equation}
	\label{eq:LinObs}
	b^{n} = AX^{n}+QV^{n},
\end{equation}
where $A$ is a real $k\times rm$ matrix. The implicit filter takes on a simple structure since, from Eq. (\ref{eq:F}), we get
\begin{equation}
	F_j(x)=\frac{1}{2}(x-\mu_j)^T\Sigma_j^{-1}(x-\mu_j)+\phi_j,
\end{equation}
with
\begin{equation}
\label{eq:LinObsSig}
	\Sigma_j^{-1} = (G(X_j^n,t^n)G(X_j^n,t^n)^T)^{-1}+A^T(QQ^T)^{-1}A,
\end{equation}
\begin{equation}
\label{eq:LinObsMu}
	\mu_j = \Sigma_j ((G(X_j^n,t^n)G(X_j^n,t^n)^T)^{-1}R(X_j^n,t^n)+A^T(QQ^T)^{-1}b^{n+1}) ,
\end{equation}
\begin{equation}
	  K_j =AG(X_j^n,t^n)G(X_j^n,t^n)^TA^T+QQ^T,
\end{equation}
\begin{equation}
\label{eq:LinObsPhi}
	    \phi_j = \frac{1}{2}(b^{n+1}-AR(X_j^n,t^n))^T K_j^{-1}(b^{n+1}-AR(X_j^n,t^n)) ,
\end{equation}
With $L_j = I$, we substitute the random map (\ref{eq:Ansatz}) into (\ref{eq:UnderdeterminedEq}) to find 
\begin{equation}
	\lambda_j=\frac{\sqrt{\rho_j}}{\sqrt{\eta_j^T\Sigma_j^{-1}\eta_j}},
\end{equation}
and the Jacobian $$J  \propto \exp(-\phi_j) (\eta_j^T\Sigma_j^{-1}\eta_j)^{-n/2}.$$ Since $\Sigma_j$ is symmetric, the values that the random variable $J$ can take on are bounded above and below. The Jacobian $J$ can vary dramatically from one sample (of $\xi$, respectively $\eta$) to another, especially if the largest and smallest eigenvalues of $\Sigma_j$ are separated by a large gap. As an example, an approximation of the pdf of $J$ for $n=2$ and $\Sigma^{-1}=\mbox{diag}(1,\mbox{ }0.5)$ is shown in figure \ref{fig:BetaPDF}. 
The pdf has two peaks, at the left and right ends of the interval over which $J$ is defined. This interval is determined by the eigenvalues of $\Sigma$ and, in this example, $J$ can can take on any value in the interval $[0.5,1]$. Choosing $L_j=I$ produces a Jacobian that can vary significantly from particle to particle. The goal however is to make the weights as uniform as possible.  

If we choose $L_j$ such that $\Sigma_j= L_j^T L_j$, we find $\lambda_j =\sqrt{\rho_j}$ and  $J=\det L_j$. This Jacobian can be expected to be roughly constant as long as the particles are reasonably
close to each other, as they are expected to be with our filter. In the special case of additive noise, i.e. $G(x^n)=G=\mbox{const}.$ in (\ref{eq:SODE}), and linear
observations available at each point in time, the Jacobian $J=\det L_j$ is constant and need not be computed. In fact, the implicit filter with random maps is, for this special case and with this choice of $L_j$, equivalent to optimal importance sampling \cite{Doucet2002, Bocquet2010,DelMoral1998}. The implicit filter is thus optimal in this case in the sense that its weights have minimum variance \cite{Doucet2002}.

In the general case, we have the Hessian evaluated at the minimum, say $H_j$, at our disposal because we use Newton's method to minimize $F_j$  and thus have:
\begin{equation}
	\label{eq:ExpandF}
 F_j(X_j) = F_j(\mu_j) + \frac{1}{2}(X_j-\mu_j)^T H_j (X_j-\mu_j) + \mbox{higher order terms}.
\end{equation}
Choosing $L_j$ so that $H_j^{-1} = L_j^TL_j$ is a good choice, especially if $F_j$ is quadratic or nearly so. This choice of $L_j$ also suggests a good initialization for the numerical computation of the parameter $\lambda$ in the random map (\ref{eq:Ansatz}). One can expect $\lambda$ to be on the order of $\sqrt{\rho}$ and one chooses $\lambda_j^0 = \sqrt{\rho_j}$. In all the examples below, the minimization in (\ref{eq:TheMin}), as well as the iterative computation of $\lambda$ converged after a few steps with this set-up. 

It is also interesting to compare the random map implementation of the implicit filter to an implementation outlined in \cite{Chorin2010b}. There, the function $F_j(X_j)$ is replaced by its quadratic approximation 
\begin{equation}
	\label{eq:ApproxF}
	F_j^0(X_j) = F_j(\mu_j) + \frac{1}{2}(X_j-\mu_j)^T H_j (X_j-\mu_j).
\end{equation} 
Instead of solving (\ref{eq:UnderdeterminedEq}), one solves $F_j^0(X_j)-\phi_j^0=0.5\xi_j^T\xi_j$, where $\phi_j^0 = \min F_j^0$ can be computed by formulas similar to (\ref{eq:LinObsSig})-(\ref{eq:LinObsPhi}). The solution of this approximate equation can be obtained by a Cholesky decomposition of $H_j$ (the Hessian of $F_j$ at the minimum) and the Cholesky decomposition also yields the Jacobian $J$. A reweighting of the particles to account for the fact that one solves an approximate equation rather than (\ref{eq:UnderdeterminedEq}) gives the weights
\begin{equation}
	\label{eq:ApproxWeights}
	w_j^{n+1}\propto w_j^{n}\exp(\phi_j) \cdot \exp(F_j(X_j)-F_j^0(X_j)) \cdot J.
\end{equation} 
The extra term $\exp(F_j(X_j)-F_j^0(X_j))$ can produce low weights if the quadratic approximation of $F^0_j$ is not close to $F_j$ in the neighborhood of the sample $X_j$. The random map (\ref{eq:Ansatz}) eliminates this factor because one solves equation (\ref{eq:UnderdeterminedEq}), rather than an approximate equation. We traded an iterative solution of a scalar equation for possibly small weights due to a quadratic approximation. 

\section{Filtering a stochastic Lorenz attractor}
\label{sec:Lorenz}
The stochastically driven Lorenz attractor \cite{Lorenz63} has been used as a testbed for data assimilation algorithms on many occasions \cite{Chorin2004,Miller1999,Miller1994}. We follow this trail and test the implicit filter on the stochastic Lorenz attractor with additive noise
\begin{eqnarray}
\label{eq:Lorenz1}
	dx &=& \sigma(y-x)dt  + g_1 dW_1,  \\
\label{eq:Lorenz2}
	dy &=& (x(\rho-z)-y)dt + g_2 dW_2,  \\
\label{eq:Lorenz3}
	dz &=& (xy-\beta z)dt  + g_3 dW_3, 
\end{eqnarray}
with the standard parameters $\sigma = 10$, $\rho=28$, $\beta=8/3$, and initial conditions $x(0)=-5.91652$, $y(0)=-5.52332$, $z(0)=24.5723$. The noise is chosen equally strong for all variables. Specifically, we choose $g_1=g_2=g_3=g=\sqrt{2}$.

\subsection{Discretization of the dynamics}
\label{sec:LorenzDiscretization}
\begin{figure}[htbp]
	\begin{center}
		{\includegraphics[width=0.8\textwidth]{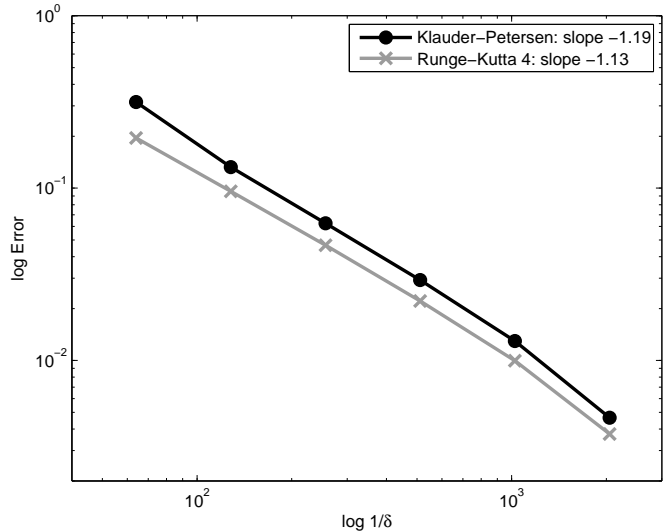}}
	\caption{Convergence of Klauder-Petersen scheme and Runge-Kutta scheme for the stochastic Lorenz attractor.}
	\label{fig:LorenzConvergence}
	\end{center}
\end{figure}

\begin{figure}[htbp]
	\begin{center}
		{\includegraphics[width=1\textwidth]{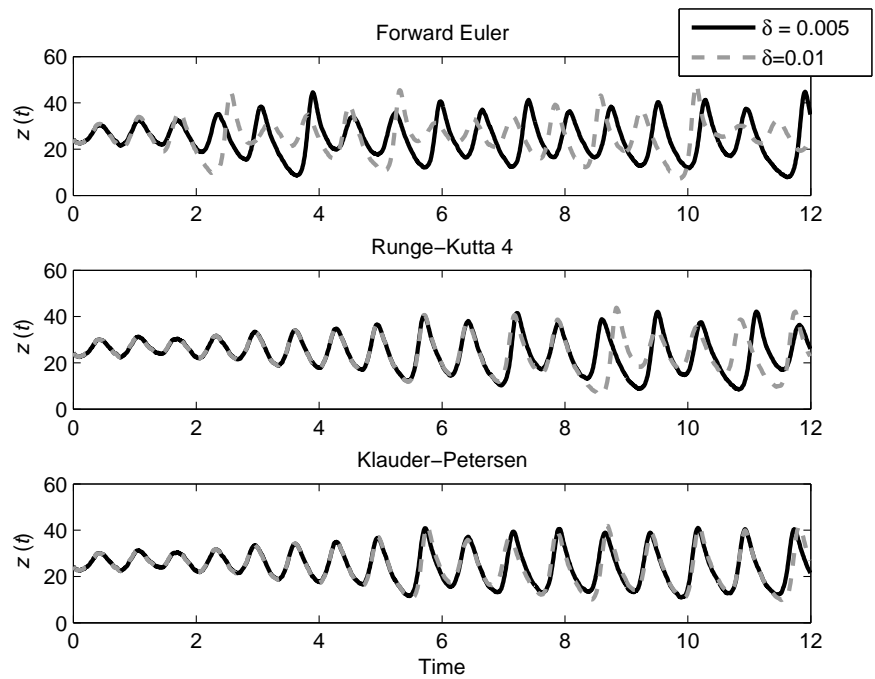}}
	\caption{Discrete time approximation of the $z$-variable for a given BM and using the Euler scheme (top), Runge-Kutta 4 (middle) and Klauder-Petersen scheme (bottom).}
	\label{fig:LorenzConvergenceTime}
	\end{center}
\end{figure}

We discretized the continuous equations by the Klauder-Petersen (KP) scheme \cite{Kla1}
\begin{eqnarray}
x^{n+1,*}&=&x^n+\delta f(x^n)+g\Delta W_1,\\
x^{n+1}&=&x^n+\frac{\delta}{2}\left(f(x^n)+f(x^{n+1,*})\right)+g\Delta W_2,
\end{eqnarray}
where $\Delta W_1, \Delta W_2 \sim \mathcal{N}(0,\delta I)$ and where $f(\cdot)$ can be read off the Lorenz attractor (\ref{eq:Lorenz1}) - (\ref{eq:Lorenz3}).
The scheme is second-order accurate for $g=0$. With $g\neq0$ (additive noise), the scheme is of strong order 1, i.e. the mean of the error at $t=t^n$ is bounded by $\mbox{Const}\cdot\delta$. Figure~(\ref{fig:LorenzConvergence}) shows the convergence of the KP scheme for the stochastic Lorenz attractor after one dimensionless time unit. 
The graph (black line) shows the mean of the error as a function of the time step. The error was approximated by the difference between the solution with time step $\delta$ and the reference solution with time step $\delta_{ref} = 2^{-12}$. The mean of the error norms (not the difference in mean error!) was approximated by running 1000 experiments and averaging. We observed the expected first-order decay in the mean of the error. For comparison, we also computed the convergence of a fourth-order Runge-Kutta (RK) scheme, where we added a Gaussian with variance $\delta g^2 I$ after each full step. For $g=0$, this scheme is fourth-order. For $g\neq0$, it is of strong order 1 (because no integrals of the BM are evaluated \cite{KloedenPlaten}). We ran 1000 experiments with $g=\sqrt{2}$ to approximate the mean of the error and observed the expected first-order stochastic convergence (light-grey line in figure \ref{fig:LorenzConvergence}). The stochastic orders of convergence for the KP or RK schemes were thus no better than that of the simple forward Euler scheme \cite{KloedenPlaten}. However, the forward Euler discretization could not follow the solution of the SDE for large integration times, because of its low accuracy in the deterministic part. The situation is illustrated in figure~\ref{fig:LorenzConvergenceTime}.
The forward Euler discretization diverged after roughly 2 dimensionless time units, while the KP and RK schemes converged for the significantly longer integration time of 12 units. We should point out that figure~\ref{fig:LorenzConvergenceTime} shows only one representative model run and that each scheme evolved under a different BM.

In data assimilation applications, one should aim at a scheme that is of low order (to make the filter efficient and fast), but at the same time catches the dynamics of the underlying SDE (the observation data may be incompatible with an under-resolved discretization). For the Lorenz attractor, the high deterministic accuracy of the RK scheme appeared unnecessary. The simpler KP scheme yielded a comparable stochastic convergence and followed the solutions of the underlying SDE for long enough. We were thus content with the KP scheme and a time step $\delta=0.01$.

\subsection{Filtering results}
\label{sec:LorenzResults}
We start by considering a case with observations of all three state variables at every time step. The observations were corrupted by noise with variance $0.1$, i.~e. we chose $Q=\sqrt{0.1}I$ in (\ref{eq:Observation}). We applied the implicit filter as explained in section \ref{sec:Construction}. The function $F_j$, as given by equation~(\ref{eq:KPF}), was minimized by Newton's method, initialized by a model run of one step without noise. 
 
At each step, we sampled a 6 dimensional standard normal variate $\xi_j$ (the reference variable) and computed the random direction $\eta_j=\xi_j \sqrt{\xi_j^T\xi_j}$ to be used in the random map. As explained in section \ref{sec:Construction}, we chose $L_j$ in (\ref{eq:Ansatz}) to be a Cholesky factor of the Hessian evaluated at the minimum. Substitution of the map~(\ref{eq:Ansatz}) into the algebraic equation (\ref{eq:UnderdeterminedEq}) gave the required equation for $\lambda_j$, which we solved by a Netwon method. The iteration was initialized by choosing $\lambda_j^0=\sqrt{\rho_j}$ and typically converged within 4-6 steps. Finally, we computed the weight of the particle using (\ref{eq:WeightMap}) and the numerical derivative $\partial \lambda / \partial \rho$, with a perturbation $\Delta \lambda = 10^{-5}\sqrt{\rho}$. We repeated this process for each particle and resample with ``algorithm~2'' in~\cite{Doucet2002}. We decided to resample at every time an observation becomes available.

We compared the implicit filter with an SIR filter \cite{Doucet2002,Doucet2001,Gordon1993}. To that effect, we ran 1000 twin experiments. That is, we ran the model for $1200$ time steps to produce artificial observations corrupted by the assumed noise. This model run was the reference we wished to reconstruct using the SIR and the implicit filters. For each experiment, the error at time $t^n$ is measured by 
\begin{equation}
 e^n = ||x^n_{ref}-x^n ||,
\end{equation}
where the norm is Euclidean, $x^n_{ref}$ is the reference state, and $x^n$ is the reconstruction by a filter. We computed this error after 5, 10 and 12 dimensionless time units (i.e. after 500, 1000 and 1200 steps) for both filters. We then computed the mean value of the error norms (mean error, for short) and the mean of the variance of the error norms (mean variance of the error, for short). The mean of the error norm is a better estimate of than the mean error because it does not allow for cancellations. The mean variance of the error is not the variance of the mean, it is a fair estimate of the error in each individual run. Our results are in table \ref{tab:LorenzLinObs}.
\begin{table}[htdp]
\caption{Results for observations of the state at every time step.}
\begin{center}
\begin{tabular}{c c c c}
\hline
 \# of Particles & \multicolumn{3}{c}{Mean error / mean variance of the error, implicit filter} \\
    &       $t=5$    &    $t=10$     & $t=12$ \\
 5  &   0.4146/0.2624 & 0.4369/0.3687 & 0.4270/0.3216 \\
 10 &   0.3215/0.1351 & 0.3289/0.1391 & 0.3311/0.1690 \\
 20 &   0.2783/0.0979 & 0.2822/0.1018 & 0.2866/0.0991\\
 30 &   0.2691/0.0914 & 0.2728/0.0931 & 0.2688/0.0908\\
    & \multicolumn{3}{c}{Mean error / mean variance of the error, SIR filter} \\
    &       $t=5$    &    $t=10$     & $t=12$ \\
 5  &  0.7915/1.8066 & 1.1751/4.0425 & 1.2544/4.3517 \\
 10 &  0.4464/0.6503 & 0.0511/0.9587 & 0.4158/0.4158 \\
 20 &  0.3159/0.1783 & 0.3196/0.2920 & 0.3156/0.1815 \\
 30 &  0.2798/0.1016 & 0.2838/0.1013 & 0.2810/0.0999 \\
 50 &  0.2695/0.0910 & 0.2688/0.0919 & 0.2711/0.0913 \\
\hline
\end{tabular}
\end{center}
\label{tab:LorenzLinObs}
\end{table}

Table \ref{tab:LorenzLinObs} shows that the implicit filter produced a small mean error and small mean variance of the error with 20-30  particles.  We can also see that the statistics converged with about 20 particles. With 20 particles the mean error variance is of the order of the variance of the observations. Even 10 particles yielded good results. The SIR filter required about 50 particles to yield comparable accuracy. For either filter, we observed a significant increase in the mean error and mean variance of the error with increasing time if the number of particles is too low (about 5 for implicit filter, about 20 for SIR filter). The increase in mean error and mean error variance was due to sample impoverishment: as time progresses, the quality of the particle ensemble decreased, i.~e. more and more particles had low weights. The effects of sample impoverishment were less severe for the implicit filter than for the SIR filter. 

Table~\ref{tab:LorenzXLinObs} shows error statistics for the SIR and implicit filters when observations of the $x$-variable only are available, i.~e. the observations are dense in time, but ``sparse in space.''  The observations were corrupted by noise with variance 0.1. The results are qualitatively the same as above. The implicit filter required about 20 particles, while the SIR filter needed about 50 particles for comparable accuracy. 
\begin{table}[htdp]
\caption{Results for observations of the variable $x^n$ at every time step.}
\begin{center}
\begin{tabular}{c c c c}
\hline
 \# of Particles & \multicolumn{3}{c}{Mean error / mean variance of the error, implicit filter} \\
    &       $t=5$    &    $t=10$     & $t=12$ \\
 5  &  0.4846/0.2624 & 0.4369/0.3687 & 0.4270/0.3216 \\
 10 &  0.3215/0.1351 & 0.3284/0.1391 & 0.3311/0.1609 \\
 20 &  0.2783/0.0979 & 0.2822/0.1018 & 0.2806/0.0991 \\
 30 &  0.2691/0.0914 & 0.2728/0.0931 & 0.2688/0.0908 \\

    & \multicolumn{3}{c}{Mean error / mean variance of the error, SIR filter} \\

    &       $t=5$    &    $t=10$     & $t=12$ \\
 5  &  0.7915/1.8066 & 1.1750/4.0425 & 1.2544/4.3517 \\
 10 &  0.4464/0.6503 & 0.5011/0.9587 & 0.4311/0.4158 \\
 20 &  0.3159/0.1783 & 0.3196/0.2920 & 0.3156/0.1815 \\
 30 &  0.2798/0.1016 & 0.2838/0.1013 & 0.2810/0.0999 \\
 50 &  0.2693/0.0910 & 0.2688/0.0919 & 0.2711/0.0913 \\
\hline
\end{tabular}
\end{center}
\label{tab:LorenzXLinObs}
\end{table}

Finally, we considered the case of observations that are sparse in time. Observations of all three state variables, corrupted by noise with variance 0.1, became available every $0.48$ dimensionless time units (every 48 steps). This is a hard data assimilation problem and some filters miss transitions from one wing of the Lorenz butterfly to the other~\cite{Miller1999}. 

The larger dimension of this problem required an additional tweak of the algorithm. The problem is of dimension 288: 3 dimensions for the Lorenz attractor, times 2 for the intermediate step $x^*$ of the KP scheme, times 48 for the gap between observations. If the variance matrix of the reference variable $\xi$ is the identity matrix $I$, we are expressing a vector variable $X_j$ of small variance as a function of a unit reference variable, and this produces very small Jacobians $J$ which can lead to underflow. One solution is to rescale $\xi$ which, after all, is arbitrary. What we did instead is keep track of the logarithms of the weights rather than the weights themselves wherever we could; this solved the problem.

\begin{figure}[h!tbp]
\label{fig:Noise}
	\begin{center}
		{\includegraphics[width=1\textwidth]{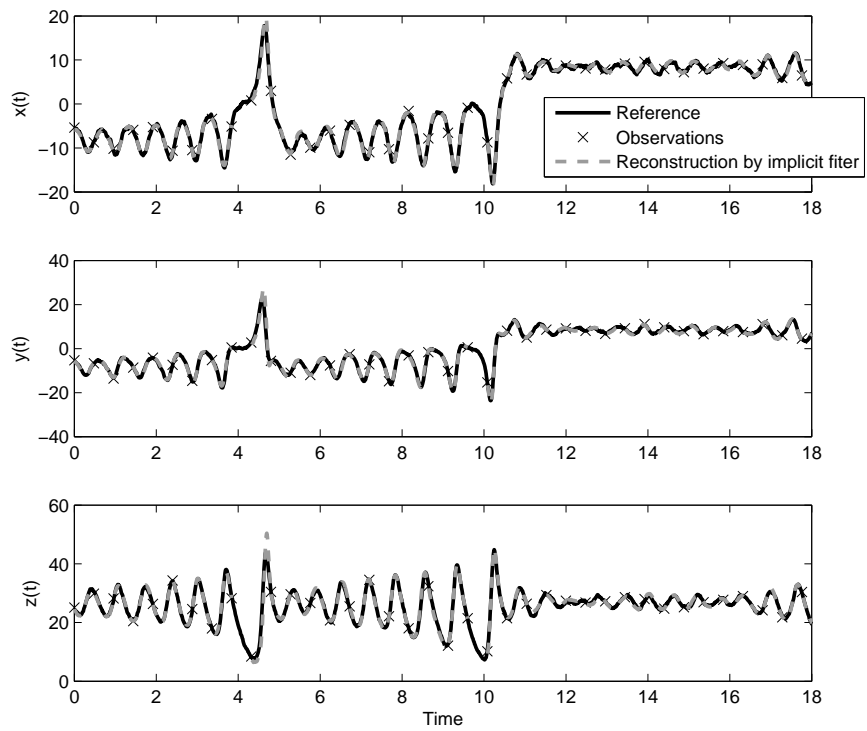}}
	\caption{A twin experiment with an implicit filter and 10 particles.}
	\label{fig:LorenzExperiment}
	\end{center}
\end{figure}

We ran 1000 twin experiments with this algorithm. Table \ref{tab:LorenzSparseObs} shows the error statistics.
\begin{table}[h!tdp]
\caption{Results for sparse observations of the state.}
\begin{center}
\begin{tabular}{c c c c}
\hline
 \# of Particles & \multicolumn{3}{c}{Mean error / mean variance of the error for implicit filter} \\
    &       $t=4.8$    &    $t=9.6$ \\
 5  &  0.1924/0.1750 & 0.2192/0.3457  \\
 10 &  0.2101/0.4103 & 0.2317/0.4905  \\
 20 &  0.1676/0.0523 & 0.1927/0.1646  \\

    & \multicolumn{2}{c}{Mean error / mean variance of the error for SIR filter} \\

    &       $t=4.8$    &    $t=9.6$    \\
 10  &  0.6508/1.0093 & 0.9964/1.9970 \\
 20  &  0.4313/0.4663 & 0.5352/0.7661 \\
 50  &  0.3368/0.2594 & 0.4271/0.5445  \\
 100 &  0.2156/0.0929 & 0.2336/0.1229 \\
\hline
\end{tabular}
\end{center}
\label{tab:LorenzSparseObs}
\end{table}
Results of one of the twin experiments are shown in figure~\ref{fig:LorenzExperiment}.
The implicit filter yielded good results with 20-30 particles, while an SIR filter required about 100 particles for  comparable accuracy (our results for the SIR filter are in agreement with those reported in \cite{Chorin2004,Miller1999}). 
We had problems with our minimization algorithm for gaps that exceed $0.5$ time units. A more sophisticated initialization or a more robust minimization can provide a cure. A detailed discussion of these issues will be taken up in future papers.

\subsection{Discussion}
The SIR and implicit filters both worked well on the stochastic Lorenz attractor and, with a sufficient number of particles, reconstructed the reference solution reliably. In all cases we considered, we observed that the implicit filter required fewer particles than the SIR filter to give a comparable accuracy. We observed that the ``focussing'' of the particles towards the observations was most beneficial when the gap between observations is large. This is indicated by the larger number of particles required in SIR than in the implicit filter. The reason is that the unguided SIR particles are very likely to become unlikely when the gap between observations is large, so that the SIR importance density and the target density can become nearly mutually singular \cite{Chorin2010b,Bickel2008,Snyder2008}. The particles of the implicit filter on the other hand are guided towards the observations because they are generated by solving (\ref{eq:UnderdeterminedEq}), which incorporates information from the available data.

The computational cost of these filters is comparable for the examples of the stochastic Lorenz attractor. The implicit filter requires fewer particles, but the computations for each particle are more expensive when compared to the SIR filter. The random map solution of the algebraic equation (\ref{eq:UnderdeterminedEq}) is efficient and reliable. 

\section{Filtering a stochastic Kuramoto-Sivashinsky equation}
\label{sec:KS}
The Kuramoto-Sivashinksy equation \cite{Kuramoto1975,Sivashinsky1977} is a chaotic partial differential equation that models laminar flames or reaction-diffusion systems (see \cite{Hu2001,Hyam1986}). Recently, stochastic Kuramoto-Sivashinsky (SKS) equations have also been used as a large dimensional test-problem for data assimilation algorithms \cite{Chorin2004,Jardak2009}. We follow in these footsteps and test the implicit filter with random maps on the SKS equation
\begin{equation}
	\label{eq:KS}	
	u_t +uu_x+u_{xx}+\nu u_{xxxx}=g\mbox{ }W(x,t)
\end{equation}
where $\nu>0$ is the viscosity, $g$ is a scalar and $W(x,t)$ is a stochastic process. We restrict ourselves to the strip $x\in[0,L], t\geq0$ and consider the case of $L$-periodic boundary conditions. Expanding the solution into a Fourier series $u(x,t) = \sum \tilde{U}_k(t) \exp(i\omega_k x)$ transforms (\ref{eq:KS}) into an infinite dimensional stochastic ordinary differential equation of the form
\begin{equation}
	\label{eq:SPDE}	
	d\tilde{U}=(\mathcal{A}(\tilde{U})+\mathcal{F}(\tilde{U}))dt+ g\mbox{ } dW_t,
\end{equation}	
where $\mathcal{A}$ is a diagonal linear operator and where $\mathcal{F}$ is a nonlinear operator. We assume that the noise process $dW_t$ is a cylindrical Brownian motion \cite{Jentzen2009}, i.e., there exists a sequence $q_n\geq 0$, $n\geq 1$, of positive real numbers, and a real number $\gamma \in (0,1)$ such that 
\begin{equation}
	\label{eq:CBMCondition}	
	\displaystyle \sum_{n=1}^{\infty} \lambda_n^{2\gamma -1}q_n < \infty,
\end{equation}	
where $\lambda_n$ are the eigenvalues of $\mathcal{A}$ in (\ref{eq:SPDE}). Let $\beta_t^n$ be independent BM's. The cylindrical BM $dW_t$ is given by the infinite series
\begin{equation}
	\label{eq:CBM}	
	dW_t =  \sum_{n=1}^{\infty}\sqrt{q_n}e_n\beta_t^n,
\end{equation}	
where the $e_n$ are unit vectors. The coefficients $q_n$ control the continuity of the noise process $dW_t$ in space. For example, $q_n=1$ for all $n$ corresponds to space-time white noise, and exponentially decaying $q_n$'s make the noise continuous in space while it remains white in time. We will consider two noise processes to drive the KS equation, namely space-time white noise and, following  \cite{Chorin2004,Chueshov2000,Shirikyan2002, LordRougemont},  spatially smooth noise with $q_n = \exp(-\left|\omega_n\right|)$. In figure~\ref{fig:Noise}, a realization of space-time white noise is shown in comparison with a realization of a spatially smooth noise process. 
\begin{figure}[htbp]
\label{fig:Noise}
	\begin{center}
		{\includegraphics[width=1\textwidth]{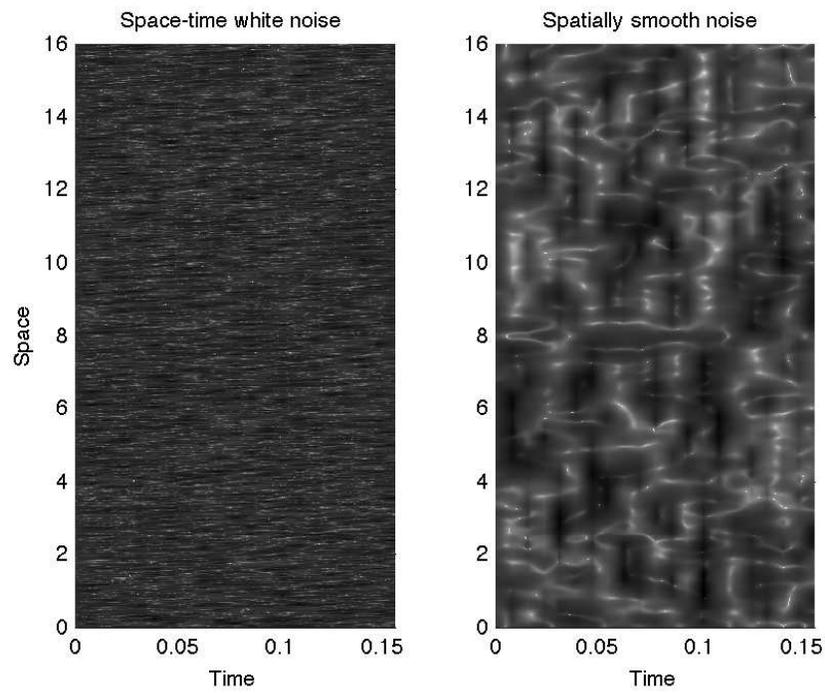}}
	\caption{A realization of space-time white noise (left) compared to a realization of spatially smooth noise (right).}
	\label{fig:Noise}
	\end{center}
\end{figure}
The absence of any correlation in space or time is evident from the left panel of figure \ref{fig:Noise}. The right panel illustrates the continuity of the noise in space at a fixed time.

A projection of equation (\ref{eq:SPDE}) onto an {\it m}-dimensional subspace spanned by {\it m} Fourier modes $\tilde{U}_k(t)$ yields an {\it m}-dimensional It\^{o}-Galerkin approximation of the SKS equation 
\begin{equation}
	\label{eq:ItoGalerkin}	
	dU=(\mathcal{L}(U)+\mathcal{N}(U))dt+ g\mbox{ }  dW_t^m,
\end{equation}	
where $U$ is a finite dimensional column vector whose components are the Fourier coefficients of the solution  and where $dW_t^m$ is a truncated cylindrical BM, obtained by projection of the cylindrical BM $dW_t$ into the Fourier modes. Assuming that the initial conditions $u(x,0)$ are odd with $\tilde{U}_0(0)= 0$ and that $dW_t^m$ is imaginary, all Fourier coefficients $U_k(t)$ are imaginary for all $t\geq0$. Writing $U_k=i\hat U_k$ and subsequently dropping the hat gives
\begin{eqnarray}
	\mathcal{L}(U) = \mbox{diag}(\omega_k^2-\nu\omega_k^4)U, \\
	\left\{\mathcal{N}(U)\right\}_k =  -\frac{\omega_k}{2}\sum_{k'=-m}^m U_{k'}U_{k-k'},
\end{eqnarray}
where $\omega_k = 2\pi k/L$, $k=1,\dots,m$ and $\left\{\mathcal{N}(U)\right\}_k$ denotes the $k^{th}$ element of the vector $\mathcal{N}(U)$.  We choose a period $L=16\pi$ and a viscosity $\nu = 0.251$, to obtain SKS equations with 31 linearly unstable modes. This set-up is similar to the SKS equation considered in \cite{Jardak2009}. With our parameter values there is no steady state as in \cite{Chorin2004}. We chose zero initial conditions $U(0) = 0$, so that the solution evolves solely due to the effects of the noise. 

\subsection{Numerical integration and convergence}
The 
numerical integration of (\ref{eq:ItoGalerkin}) is more delicate than for the Lorenz attractor.  The forward Euler scheme is unstable for any reasonable time step so that one must consider more sophisticated schemes to discretize (\ref{eq:ItoGalerkin}), see e.~g. \cite{Jentzen2009,KloedenPlaten,LordRougemont, Milshtein1978,Milstein1998}. We found that fully implicit schemes, for example implicit Euler or an implicit 1.5-strong-order scheme, are numerically awkward for the SKS equation (and, in fact, for most high dimensional problems). The exponential Euler scheme \cite{Jentzen2009} can be thought of as a stochastic version of exponential time differencing \cite{Cox2002} and is tailor-made for nonlinear equations whose stiffness arises from their linear parts. While the scheme is only first order in time, competing schemes, for example the linear-implicit Euler or the Lord-Rougemont  \cite{LordRougemont} schemes, converge even slower. Similar observations were made in \cite{Jentzen2009}. Taking into account both the time discretization and space truncation error, the exponential Euler scheme appeared superior to other schemes we considered. For the SKS equation, this scheme takes the form
\begin{equation}
	\label{eq:Scheme}
	U^{n+1} = e^{B\delta}U^n+B^{-1}(e^{B\delta}-I)\mathcal{N}(U^n)+g\sqrt{0.5DB^{-1}(e^{2B\delta}-I)}\Delta W^{n+1}
\end{equation}
where $B=\mbox{diag}(\omega_k^2-\nu\omega_k^4)$ and $D=\mbox{diag}(q_k)$. Note that all the matrices are diagonal, so that the numerical integration can be implemented efficiently, even if $m$ is large. In the notations of Eq.~(\ref{eq:SODEDiscrete}), we write
\begin{eqnarray}
	R(U^n,t^n) &=& = e^{B\delta}U^n+B^{-1}(e^{B\delta}-I)\mathcal{N}(U^n), \\
	G(U^n) &=&  g \sqrt{0.5QB^{-1}(e^{2B\delta}-I)}. \label{eq:NoiseVariance}
\end{eqnarray}

\begin{figure}[htbp]
	\begin{center}
		{\includegraphics[width=1.0\textwidth]{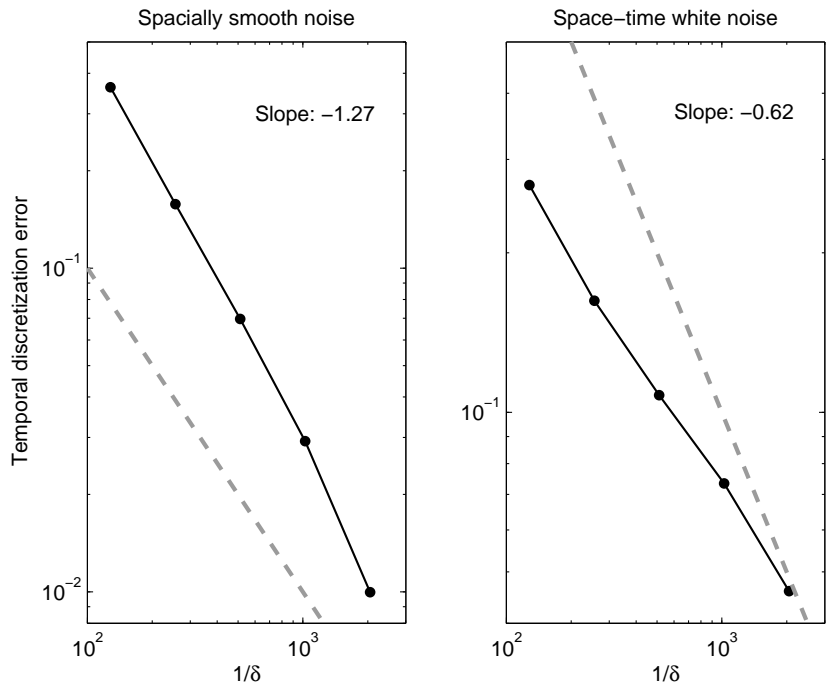}}
	\caption{Convergence of the exponential Euler scheme. Left: spatially smooth noise. Right: space-time white noise. Solid line: mean error at $T=20$. Broken line: order line 1.0.}
	\label{fig:Convergence}
	\end{center}
\end{figure}

\begin{figure}[htbp]
	\begin{center}
		{\includegraphics[width=1.0\textwidth]{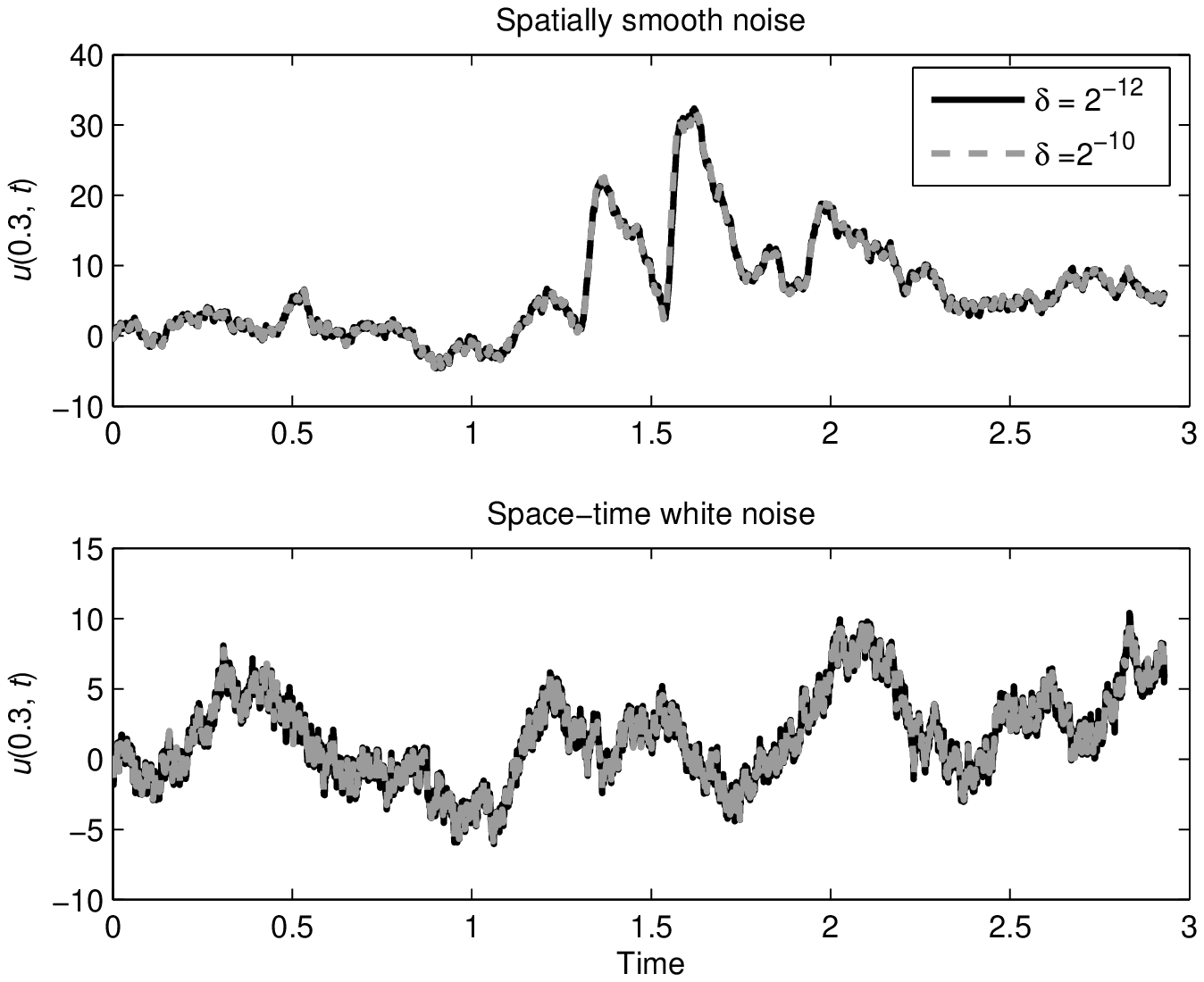}}
	\caption{Discrete approximation of $u(\hat{x}=0.3,t)$ for different noise processes and different time steps. For each case, the numerical approximations share the same realization of the noise.}
	\label{fig:TimeConvergence}
	\end{center}
\end{figure}

To assess the convergence of the exponential Euler scheme we calculated a very accurate reference solution with a time step of $2^{-12}$ and compared it to approximations with varying time steps. The number of Fourier modes was held fixed: 512 in the case of spatially smooth noise and 1024 in the case of space-time white noise. The mean error was approximated as the average of $||\hat{u}(x,T)-u_{ref}(x,T)||$ over 2000 experiments, where $||\cdot||$ was the Euclidean norm, $T=3$, $g=4$, and $u_{ref}$ was the reference solution. Figure~(\ref{fig:Convergence}) shows the results.
For spatially smooth noise we observed a convergence rate of about one, as expected. The scheme converged slower when we made the noise white in space, i.~e. increased the noise in high frequency modes. The exponential Euler scheme converged when the noise is white in space and time because the elements of the diagonal matrix multiplying the BM in (\ref{eq:SODE}) became smaller as the number of Fourier modes increases (see equation (\ref{eq:NoiseVariance})). Figure \ref{fig:TimeConvergence} shows the results of one of our experiments and indicates that the discretization follows the solution of the SPDE long enough for our purposes.

We were content with a time step $\delta = 2^{-10}$ and $m=128$ modes for spatially smooth noise, and $\delta=2^{-10}$, $m=512$ for space-time white noise.

\subsection{The observations}
We are solving the SKS equations in Fourier variables, but we choose to observe in physical
space (as is maybe physically reasonable). The solution of the algebraic equation
(\ref{eq:UnderdeterminedEq}) is easiest when the function $F$ is nearly diagonal, i.e., when
its linearizations around a current state are nearly diagonal matrices; this requires in particular 
that the variables that are observed coincide with the variables that are evolved by the
SDE. Observing in physical space while computing
in Fourier space creates the opposite situation, in which each observation is related to
the variables one computes by a dense matrix. Our solution of this problem demonstrates 
the effectiveness of the random map algorithm.

Specifically, 
we collected observations $h(u(x,t))$, corrupted by noise with unit variance, at the discrete locations $x_1,\dots,x_{m/2}$. Equation (\ref{eq:Observation}) becomes
\begin{equation}
	b^n=h((u(x_1,t^n),\dots,u(x_{m/2},t))^T = h(-2\mbox{Im}(E U^n))+V^n,
\end{equation} 
where $V^n\sim \mathcal{N}(0,I)$ and where $E$ is a $(2n+1)\times m$ matrix with rows $$E_j = (e^{i\omega_m x_j},e^{i\omega_{m-1} x_j},\dots,e^{i\omega_0x_j}).$$ 
For simplicity, we chose to collect the data at $m/2$ equidistant locations.

\subsection{Numerical results}
\begin{figure}[htbp]
	\begin{center}
		{\includegraphics[width=1.0\textwidth]{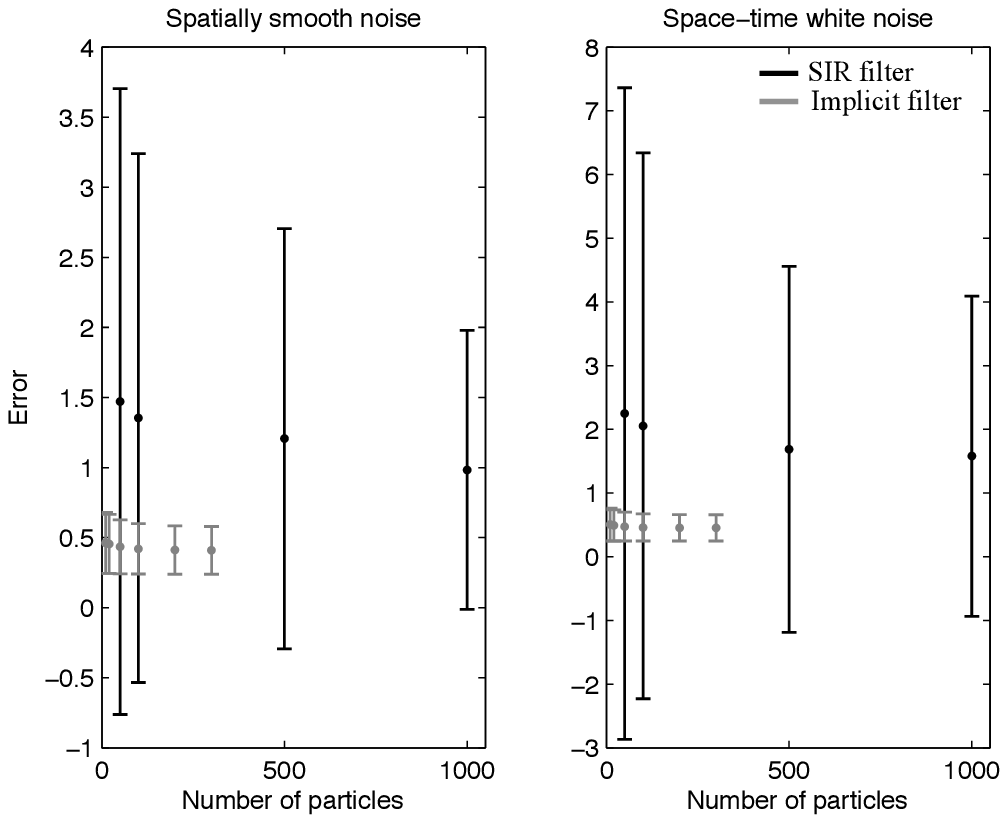}}
	\caption{Illustration of filtering results for a linear observation operator: the error statistics are shown as a function of the number of particles for SIR (black) and implicit filter (light grey).}
	\label{fig:BarLinObs}
	\end{center}
\end{figure}

To test the implicit filter we ran twin experiments as in section~\ref{sec:LorenzResults}. The error at time $t^n$ is defined as
\begin{equation}
\label{eq:KSError}
e^n=||U_{ref}^n-U_F^n||
\end{equation}
where the norm is the Euclidean norm $\left|\right| x \left|\right|=\sqrt{x^Tx}$; $U_{ref}^n$ denotes the set of Fourier coefficients of the reference run and $U_F^n$ denotes the reconstruction by the filter, both at the fixed time $t^n$. 
Table \ref{tab:KSLinObs} shows the results of 500 twin experiments for $n=100$, $g=4$ and with linear observations $h(x_j^n)=x_j^n$ at every step. The results are graphically summarized in figure \ref{fig:BarLinObs}.  Since $G$ in (\ref{eq:SODEDiscrete}) and $Q$ in (\ref{eq:Observation}) are independent of the state or time, the Cholesky factorization could be done off-line, i.~e. needed to be computed only once. It follows that $\det L$ in (\ref{eq:WeightMap}) needed not to be computed, since it is the same for all particles.
\begin{table}[htdp]
\caption{Results for linear observations at every step.}
\begin{center}
\begin{tabular}{c c c}
\hline
 & \multicolumn{2}{c}{Spatially smooth noise} \\
 \# of Particles      &   Implicit filter    &    SIR filter 	\\
 10   &  0.462345/0.217435 		& -/- 		  	\\
 20   &  0.455133/0.210594		& -/-  			\\
 50   &  0.434861/0.192192 		& 1.47129/2.23284  \\
 100  &  0.420063/0.179344 		& 1.35330/1.88725  \\
 200  &  0.41221/0.1725600 		& -/-  			\\
 300  &  0.40919/0.1700570 		& -/-  			\\
 500  &  -/-           		& 1.20573/1.498450 \\
 1000 &  -/-			  		& 0.98354/0.995908 \\
 \multicolumn{3}{c}{  }\\ 
 			& \multicolumn{2}{c}{Space-time white noise}\\
 \# of Particles&   Implicit filter    	&    SIR filter \\
10   &  0.505932/0.258586 	& -/-  \\
20   &  0.491701/0.244227	 	& -/- \\
50   &  0.473583/0.225954		& 2.24747/5.11504 \\
100  &  0.460124/0.212956 	& 2.05649/4.28421  \\
200  &  0.455131/0.208389 	& -/-  \\
300  &  0.452730/0.205857 	& -/-  \\
500  &  -/-           		& 1.68514/2.87233 \\
1000 &  -/-			  		& 1.57808/2.51565  \\
\hline
\end{tabular}
\end{center}
\label{tab:KSLinObs}
\end{table}

From table \ref{tab:KSLinObs} and figure \ref{fig:BarLinObs} we observe that the implicit filter gave very accurate results with only ten particles. The error statistics had converged,  so that it was unnecessary to perform experiments with more than 300 particles. The SIR filter collapsed (all weights were zero up to machine precision) unless the number of particles was greater than or equal to 50. Experiments with 50, 100, 500 and 1000 particles showed that SIR filter could not yield an accuracy close to that of the implicit filter with 10 particles. Even with 1000 particles the SIR filter yielded four times the mean variance of the implicit filter with ten particles. 

One can check that with our parameter choices the ratio of model-to-observation noise is larger for space-time white noise than for spatially smooth noise. The higher level of noise creates more of a problem for the SIR filter than for the implicit filter. We observed that the error of the SIR filter increased when the SKS equation was driven by white noise, while the implicit filter appeared insensitive to the nature of the model noise. In our experience, the implicit filter performs well with a large model-to-observation noise ratio.
  
\begin{figure}[htbp]
	\begin{center}
		{\includegraphics[width=1.0\textwidth]{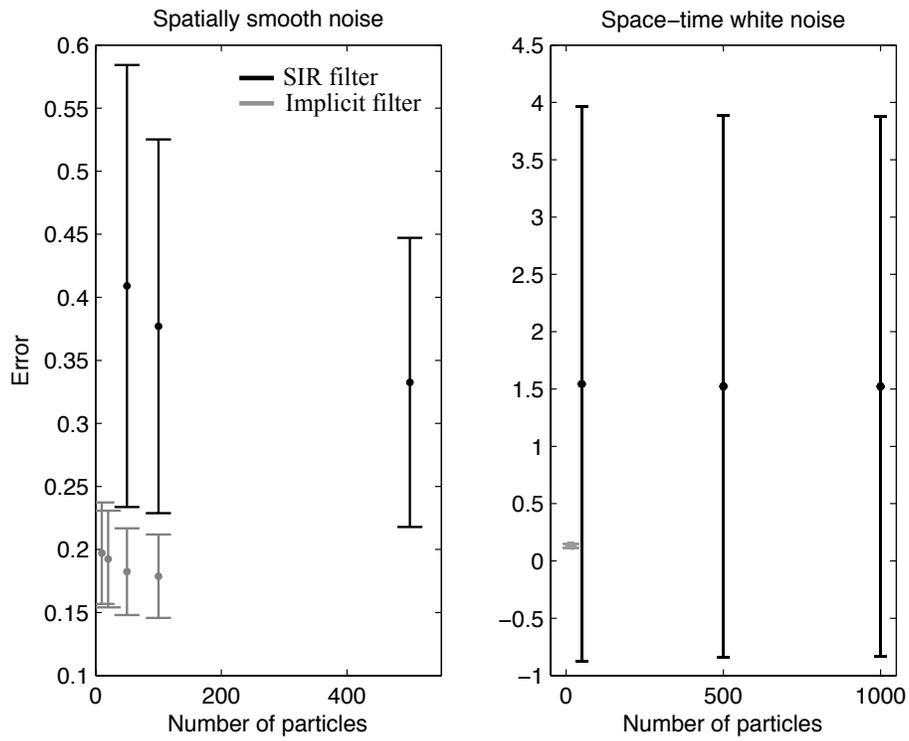}}
	\caption{Illustration of filtering results for a nonlinear observation operator: the error statistics are shown as a function of the number of particles for SIR (black) and implicit filter (light grey).}
	\label{fig:BarNonLinObs}
	\end{center}
\end{figure}

Next, we consider the nonlinear observation operator $h(x)=x+x^3$. As in section \ref{sec:LorenzResults}, the minimization of $F_j$ was done using a model run without noise as the initial guess, followed by a few full Newton steps.  Results of 500 twin experiments are shown in table \ref{tab:KSNonLinObs} and figure \ref{fig:BarNonLinObs}.
\begin{table}[htdp]
\caption{Results for nonlinear observations at every step.}
\begin{center}
\begin{tabular}{c c c c c}
\hline
& \multicolumn{2}{c}{Spatially smooth noise} \\
  \# of Particles      &   Implicit filter    &    SIR filter\\
 10   &  0.197085/0.0401874 		& -/-  			  \\ 
 20   &  0.192486/0.0383204 		& -/-  			  \\
 50   &  0.182398/0.0343374 		& 0.408985/0.175277   \\
 100  &  0.178808/0.033115		& 0.377034/0.148200    \\
 500  &  -/-           			& 0.332515/0.114040 \\
 5000 &  -/-			  			& 0.280989/0.082068\\
 \multicolumn{3}{c}{ }\\
&\multicolumn{2}{c}{Space-time white noise}\\
\# of Particles  &   Implicit filter    &    SIR filter \\ 
 10    &  0.133155/0.0181577 	& -/-  \\
 20    &  0.132795/0.0180349 & -/- \\
  50   &  -/- 					& 1.54282/2.41919 \\
 100   &  -/- 					& 2.05649/4.28421  \\
 500   &  -/-           			& 1.52291/2.36136 \\
 5000  &  -/-			  		& 1.52078/2.35841 \\
\hline
\end{tabular}
\end{center}
\label{tab:KSNonLinObs}
\end{table}
We observe from table \ref{tab:KSNonLinObs} and figure \ref{fig:BarNonLinObs} that the implicit filter outperformed the SIR filter. A SIR filter with 5000 particles gave less accurate results than the implicit filter with ten particles for either noise process. The results are similar to those obtained for a linear observation operator.

\begin{figure}[htbp]
	\begin{center}
		{\includegraphics[width=1.0\textwidth]{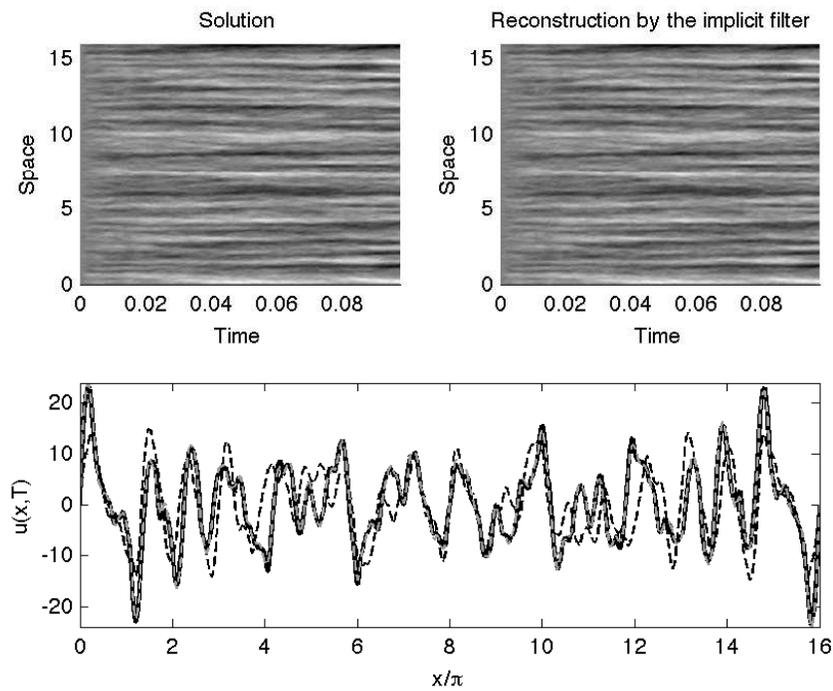}}
	\caption{Outcome of a twin experiment with data sparse in space and time. Solid black line (almost hidden): reference. Dashed gray line: implicit filter with 10 particles. Dashed black line: SIR filter with 1000 particles.}
	\label{fig:ResultsSparse}
	\end{center}
\end{figure}

Last, we consider the case of linear observations at every other time step. We ran 500 twin experiments. In each experiment we integrated the SKS equation driven by smooth noise with $g=1$ until $t^n=100\delta$. We averaged the results to estimate the error statistics. The results are shown in table \ref{tab:KSSparseObs} and figure~\ref{fig:ResultsSparse} shows of one of the twin experiments.
\begin{table}[htdp]
\caption{Results for sparse observations.}
\begin{center}
\begin{tabular}{c c c}
\hline
    & \multicolumn{2}{c}{Spatially smooth noise} \\
  \# of Particles      &   Implicit filter & SIR\\
 10    &  0.391023/0.156071 & -/-		\\ 
 20    &  0.384932/0.151217 & -/-		\\
500   &  -/- &  0.280533/0.080205		\\
1000   &  -/- & 0.271989/0.075466		\\

\hline
\end{tabular}
\end{center}
\label{tab:KSSparseObs}
\end{table}

From table \ref{tab:KSSparseObs}, we observe that the implicit filter appeared insensitive to the fact that observations were not always available. The error statistics had converged. The SIR filter required at least 500 particles to achieve a similar accuracy and often collapsed with fewer particles so that a reliable estimation of the error statistics was infeasible. The error decreased compared to table~\ref{tab:KSLinObs} because we decreased the noise in the model by setting $g=1$, rather than $g=4$.

Finally we want to comment on how our results compare to those reported in \cite{Chorin2004, Jardak2009}. In \cite{Chorin2004}, Chorin and Krause considered a SKS equation with two linearly unstable modes and successfully applied a dimensional reduction to their SIR filter. In the present paper we chose a viscosity and period that yield 31 unstable modes. Thus, the SKS equation in \cite{Chorin2004} was a lot ``nicer'' and assimilating data was easier. Jardak et al. \cite{Jardak2009} considered data assimilation for a SKS equation with 32 linearly unstable modes, however their noise is milder and the numerical integration is carried out differently. They compared the performance of an ensemble Kalman filter (EnKF) to that of an SIR filter and Maximum Likelihood Ensemble Filter methods (MLEF). Only sparse observations of the Fourier coefficients were considered. The conclusion was that EnKF outperforms SIR and MLEF for linear observations but has major drawbacks for nonlinear observation operators. For nonlinear observations, SIR gave the best results. When compared to ours, the SIR particle ensembles in \cite{Jardak2009} were smaller because of the lower noise levels and because the Fourier coefficients rather than the physical solution were observed. Nonetheless, the number of SIR particles is 70-250 and thus larger than the 10-50 particles we require for the implicit filter. A more detailed comparison of our results to those in \cite{Jardak2009} is not possible because of the different assumptions.  

\section{Conclusions}
Implicit filtering is a sequential Monte Carlo technique for nonlinear, non-Gaussian data assimilation. The implicit filter is designed to keep the number of particles required manageable by focussing the particles towards the high-probability regions.
We have presented a new implementation of an implicit particle filter in which the underdetermined algebraic equation characteristic of implicit sampling is solved efficiently via a random map. The use of the random map reduces a data assimilation problem of arbitrary dimension to a sequence of minimizations of explicitly known functions, followed by  solutions of algebraic equations.

We applied the filter in our new implementation to two challenging test problems where it performed well in comparison with a standard filter. As expected, our filter became more economical, compared with alternatives, when the dimension of the problem increased or the model noise grew. The various numerical issues
that arise as the problem size increases even further will be discussed in the context of specific applications. 

\section*{Acknowledgments}
We would like to thank our collaborators at Oregon State University, Professors Robert Miller and Yvette Spitz and Doctor Brad Weir, for helpful discussion and comments. This work was supported in part by the Director, Office of Science, Computational and Technology Research, U.S. Department of Energy under Contract No. DE-AC02-05CH11231, and by the National Science Foundation under grants DMS-0705910 and OCE-0934298.

\bibliographystyle{plain}
\bibliography{References}

\end{document}